\documentclass[12pt,twoside]{article}
\usepackage{amsmath,amsbsy,amsfonts}
\usepackage{color}
\newtheorem{lemma}{Lemma}[section]
\newtheorem{prop}{Proposition}[section]
\newtheorem{theo}{Theorem}[section]

\newtheorem{coro}{Corollary}[section]
\newtheorem{rem}{Remark}[section]
\newtheorem{claim}{Claim}[section]

\let\Section=\section
\def\section{\setcounter{equation}{0}\Section}
\begin{document}
\title{Multiple solutions for a inclusion quasilinear problem with non-homogeneous boundary \\ condition through
Orlicz Sobolev spaces  }

\author{Rodrigo C. M. Nemer and Jefferson A. Santos\thanks{ Jefferson A. Santos Partially supported by CNPq-Brazil grant Casadinho/Procad 552.464/2011-2},   }

\date{}

\pretolerance10000

\maketitle

\begin{abstract}
{\it In this work we study multiplicity of nontrivial solution for the following class of differential inclusion problems with non-homogeneous Neumann condition through Orlicz-Sobolev spaces,
$$
\left \{ \begin{array}{l}
-div\big(\phi(|\nabla u|)\nabla u\big)+\phi(|u|)u\in  \lambda\partial F(u)  \ \mbox{in}\ \Omega,\\
\frac{\partial u}{\partial \nu}\in \mu\partial G(u) \text{ on }\partial \Omega,                                    
\end{array}
\right.
$$
where $\Omega\subset{\mathbb{R}}^N$ is a domain, $N\geq 2$ and $\partial F(u)$ is the generalized gradient of $F(u)$. The main tools used are Variational Methods for Locally Lipschitz Functional and Critical Point Theory. } \\
\noindent \textit{\bf 2000 AMS Subject Classification:} 35A15,
35J25, 34A36.

\noindent \textit{\bf Key words and phrases:} quasilinear equations,
non-smooth functionals, Orlicz-Sobolev.

\end{abstract}

\section{Introduction and prerequisites}
Let $\Omega\subset{\mathbb{R}}^N$, $N\geq 2$, be a bounded domain with smooth boundary $\partial \Omega$ and consider a continuous function $\phi:(0,+\infty)\rightarrow(0,+\infty)$. For $\lambda,\mu>0$, we study existence of nonnegative solutions for the differential inclusion problem with non-homogeneous Neumann condition
$$
\left \{ \begin{array}{l}
-div\big(\phi(\mid\nabla u\mid)\nabla u\big)+\phi(|u|)u\in  \lambda\partial F(u)  \ \mbox{in}\ \Omega,\\
\frac{\partial u}{\partial \nu}\in \mu\partial G(u(x)) \text{ on } \partial \Omega,                                    
\end{array}
\right.\leqno{(P_{\lambda,\mu})}
$$
where $F, G: {\mathbb{R}}\rightarrow{\mathbb{R}}~\mbox{are locally Lipschitz}$ and
$$
\partial F(t) = \{s \in  {\mathbb{R}};~  F^{o}(t;r) \geq  s r,~r \in {\mathbb{R}}  \},
$$

where $F^{o}(t;r)$  denotes the generalized directional derivative of $t \mapsto  F(t)$ in the direction of $r$, that is
$$
F^{o}(t;r) = \limsup_{y \to t,\ s \to 0} \frac{F( y+sr) - F(y)}{s}.
$$
Analogously, we define $\partial G(t)$ and $G^{o}(t;r)$.

It is well know (see e.g. \cite{Chang, papa-1}) that if $F$ is of class $C^1$, then
$$
\partial F(t) = \{F'(t)\}.
$$
In this case, one has an equation in $(P_{\lambda,\mu})$, instead of an inclusion. We refer the reader to  Clarke \cite{clarke},  Chang \cite{Chang}, Carl, Le  and Motreanu \cite{carl} and their references. A few remarks are in order.
\vskip.2cm

\begin{rem} \label{Remark 2}  Krist\'aly, Marzantowicz, and Varga \cite{Kristaly} studied the problem $(P_{\lambda,\mu})$ for $\phi(t)=|t|^{p-2}$. By using a result of Ricceri \cite{Ricceri}, they guaranteed the existence of three critical points for a nonsmooth functional associated to the problem
$$
\left \{ \begin{array}{l}
-\Delta_p u+|u|^{p-2}u\in  \lambda\partial F(u)  \ \mbox{in}\ \Omega,\\
\frac{\partial u}{\partial \nu}\in \mu\partial G(u) \text{ on } \partial \Omega. 
\end{array}
\right.
$$
\end{rem}

The study of nonlinear partial differential equations with discontinuous nonlinearities is motivated by various real-life phenomena coming from Mechanics and Mathematical Physics. Problems from the latter had been treated by several authors, such as Pucci and Serrin \cite{Pucci}, Ricceri \cite{Ricceri,Ricceri2}, Marano and Motreanu \cite{Marano,Marano2}, Arcoya and Carmona \cite{Arcoya}, Bonanno \cite{Bn1,Bn2}, Gasi\' nski and Papageorgiou \cite{gasinski}, Krist\'aly \cite{Kristaly2}, Bonanno and Candito \cite{Bn3}, Alves and Nascimento \cite{Rubia}.
More recently, Alves, Gon\c calves and Santos \cite{Abrantes} established existence  of nontrivial solutions  for the problem
$$
\displaystyle  -\text{div} ( \phi(|\nabla u|) \nabla u) - b(u)u \in \lambda \partial F(x,u) \  \text{in} \
  \Omega,
$$
where $\lambda >0$ is a parameter and $\phi:[0,+\infty) \to [0,+\infty)$ is a $C^1$-function satisfying
$$
\leqno{(\phi_1)} \ \            \displaystyle  \lim_{s\rightarrow 0^+}s\phi (s) =0 \,\,\, \mbox{and} \,\,\, \displaystyle  \lim_{s\rightarrow +\infty} s\phi (s) = + \infty,
$$
$$
\leqno{(\phi_2)} \ \ \  s \mapsto s\phi (s)~\mbox{ is increasing in}~ [0, \infty),
$$
$$
\leqno{(\phi_3)} \ \  \   \ell \leq \frac{\phi(t)t^2}{\Phi(t)} \leq m, \   t>0, \ l,m>0 \mbox{ and } \Phi(t) = \int_0^{|t|} s \phi(s) ds,  \\
$$
where $b$ is a continous function and $F$, locally Lipschitz.

In general, $\phi$ is not a power function, not even homogeneous, but it is convex. Thus, in general, Orlicz-Sobolev spaces, rather than Sobolev spaces, are used in the study of problems like $(P_{\lambda,\mu})$. For instance, in Elasticity and Geometry, some authors (cf. Fukagai and Narukawa \cite{narukawa-2}, Dacorogna \cite{dacorogna}) deal the with de problem where $\phi$ is given by
\begin{equation} \label{(1.5)}
 \phi(t)=2\alpha (1+t^2)^{\alpha-1},~ t> 0,
\end{equation}
$1 < \alpha < \frac{N}{N-2}$. In this case the corresponding Orlicz-Sobolev space is actually equal to a Sobolev space. On the other hand, the function
$$
\phi(t)=pt^{p-2}\ln(1+t)+\frac{t^{p-1}}{t+1},~ t>0.
$$
where $\frac{\sqrt{1+4N}-1}{2} < p < N - 1$, 
gives an example where the Orlicz-Sobolev space is not equal to a Sobolev space. These remarks follow by applying a result in \cite{Rao1}.


The main result of this paper is the following

\begin{theo}\label{theorem1}
Let $F,G:\mathbb{R}\rightarrow\mathbb{R}$ be locally Lipschitz functions satisfying the conditions
\begin{description}
\item{$(F_1)$} there is $c_1>0$ such that
$$
|\xi|\leq c_1(1+b(|t|)|t|), \ \xi \in \partial F(t), \ t\in \mathbb{R},
$$
with $b:(0,+\infty]\rightarrow \mathbb{R}$ a $C^{1}$ function verifying
$$
m< b_0\leq \frac{b(t)t^2}{B(t)}\leq b_1 <l^*,\eqno{(b_1)}
$$
for all $t>0$, with
$$
(b(t)t)'>0, \ t>0, \eqno{(b_2)}
$$
and $B(t)=\int_0^{|t|}b(s)s~ds;$
\item[$(F_2)$] $\displaystyle \lim_{t\to 0}\frac{\max\{|\xi|;\xi\in \partial F(t)\}}{\phi(|t|)|t|}=0$,

\item[$(F_3)$]$\displaystyle \limsup_{|t|\to+\infty}\frac{F(t)}{\Phi(t)}\leq 0;$

\item[$(F_4)$] assume that $F(0)=0$ and there is $t_0\in\mathbb{R}\setminus\{0\}$ such that $F(t_0)>0$;

\item[$(G_1)$] there is $c_2>0$ such that
$$
|\xi|\leq c_2(1+\overline{b}(|s|)|s|), \ \xi \in \partial G(s), \ s\in \mathbb{R},
$$
where $\overline{b}:(0,+\infty)\rightarrow(0,+\infty)$ satisfies $(b_1)-(b_2)$, with \break $\overline{B}(t)=\int_0^{|t|}\overline{b}(s)s~ds,$ $b_0=\overline{b}_0$ and $b_1=\overline{b}_1$, $1 <\overline{b}_0\leq \overline{b}_1<\overline{l}^*=\frac{l(N-1)}{N-l}$.
\end{description}
Then there exists a nondegenerate compact interval $[a,b]\subset (0,+\infty)$ and a number $r>0$ such that for every $\lambda\in [a,b]$, there is $\mu_0\in (0,\lambda+1]$ such that for each $\mu\in [0,\mu_0]$, the problem $(P_{\lambda,\mu})$ has at least three distinct solutions with $W^{1,\Phi}$-norms less than $r$.
\end{theo}

Here, a function $u\in W^{1,\Phi}(\Omega)$ is a solution of the problem $(P_{\lambda,\mu})$, if there are $\xi_F\in L_{\tilde B}(\Omega)$ and $\xi_G\in L_{\tilde{ \overline B}}(\partial \Omega)$,  such that for all $v\in W^{1,\Phi}(\Omega)$ we have
$$
\int_\Omega \phi(\mid\nabla u\mid)\nabla u\nabla v +\phi(|u|)uvdx=\lambda\int_\Omega \xi_Fvdx +\mu \int_{\partial \Omega}\xi_Gvdx.
$$
Moreover,
$$
\xi_{F}(x)\in\partial F(u(x)) \text{ and } \xi_{G}(x)\in\partial G(u(x)) \text{ for a.e. }x\in \Omega.
$$

The proof of Theorem \ref{theorem1} rests upon variational techniques and consists in finding  critical points of the energy functional associated to $(P_{\lambda,\mu})$, defined for each $u\in W^{1,\Phi}(\Omega)$ as
\begin{equation}\label{energy_fuctional}
J_{\lambda,\mu}(u)=\int_\Omega\Phi(|\nabla u|)+\Phi(|u|)dx-\lambda\int_\Omega
F(u)dx-\mu\int_{\partial\Omega} G(u)dx.
\end{equation}
 Due to the generality of problem $(P_{\lambda,\mu})$,
$J_{\lambda,\mu}$ is  only locally Lipschitz continuous. The solutions are obtained applying an extension of the Three Critical Points Theorem, due to Krist\'aly, Marzantowicz and Varga \cite{Kristaly}, for locally Lipschitz functions. To apply this theorem, it was necessary, among other things, to develop a new Chain Rule Theorem and a compact embedding between $W^{1, \Phi}(\Omega)$ and Orlicz spaces defined over $\partial \Omega$. We believe these tools are new, since we could find nothing related in the literature.

This paper is divided in three more sections. In the next, we give the basics of nonsmooth functionals on Banach spaces and new Chain Rule Theorem (Theorem \ref{RegradaCadeia}); in the third section, we present nonsmooth Ricceri's theorem and some aspects of Orlicz-Sobolev spaces; finally, in Section 4, we present an application of the results shown in the previous sections to a differential inclusion problem with nonhomogeneous boundary condition and a new compact embedding (Theorem \ref{imersao1}).


\section{Basics on Nonsmooth Functionals on Banach Spaces}\label{Sec2}

Let $X$ be a Banach space and $I:X\rightarrow{\mathbb{R}}$ be a
locally Lipschitz functional, or $I\in Lip_{loc}(X,\mathbb{R})$ for
short, that is, given $u\in X$, there are an open neighborhood
$V=V_u\subset X$ and a constant $K=K_V>0$ such that
$$
\mid I(v_2) -I(v_1)\mid\leq K\Vert v_2-v_1\Vert, \ v_1, v_2\in V.
$$
 The directional derivative of $I$ at $u$ in the direction of
$v\in X$ is defined by
$$
I^0(u;v)=\displaystyle
\limsup_{h\rightarrow0,~\lambda\rightarrow0^+}\frac{I(u+h+\lambda
v)-I(u+h)}{\lambda}.
$$
 It follows that  $I^0(u;.)$ is subadditive and  positively
homogeneous, that is,
$$
I^0(u;v_1+v_2)\leq I^0(u;v_1)+I^0(u;v_2)
$$
and
$$
I^0(u;\lambda v)=\lambda I^0(u;v),
$$
where  $u,v,v_1,v_2\in X$ and $\lambda>0$. As a byproduct,
$$
\mid I^0(u;v_1)-I^0(u;v_2)\mid\leq I^0(u;v_1-v_2) \leq K\| v_1-v_2\|_X,
$$
for some  $K=K_u>0$. In addition,  $I^0(u;.)$ is continuous and convex.
 The generalized gradient of $I$ at $u$ is the set
$$
\partial I(u)=\big\{\mu\in X^*;~ \langle \mu,v\rangle\leq I^0(u;v), \ v\in
X\big \}.
$$
Since $I^0(u;0)=0$, $\partial I(u)$ is the subdifferential of
$I^0(u;0)$. Further properties, definitions and remarks are given
below (cf. \cite{Chang}, \cite{carl} for proofs):

\begin{description}
 \item{$(S_1)$}~  {$~\partial I(u)\subset X^*~  \mbox{\it is convex, nonempty and
 weak*-compact}$};

\item{$(S_2)$}~  {$~ m_I(u):=\min\big\{\|\mu\|_{X^*};~\mu \in \partial
I(u)\big\}$};

\item{$(S_3)$}~ {$\partial I(u)=\big\{I'(u)\big\},~  \mbox{if}~ I \in
C^1(X,{\mathbb{R}})$};

\item{$(S_4)$}~  {$ \mbox{\it A point}~  u_0\in X~  \mbox{\it is a critical point of}~  I \ \mbox{\it if}~\ 0\in \partial I(u_0)$};

\item{$(S_5)$}~  {$  c\in {\mathbb{R}}~ \mbox{\it is a critical value of}~ I~ \mbox{\it if there is a critical point } u_0~ \mbox{\it of}~ I~\\
~~~~~\mbox{\it such that}~ I(u_0)=c$, \ {\it and we set}}
$$
K_c=\{u\in X;\ 0\in \partial I(u),\ I(u)=c\};
$$
\item{$(S_6)$}~  {$ \mbox{\it If}~ u_0~ \mbox{\it is a local minimum of}~ I,~\mbox{\it then it is a critical point of}~ I$.}

The support function of a nonempty subset $\Sigma$ of $X^*$, $X$ being such
that $X\subset X^{**}$, is the function $\sigma_{\Sigma}:
X\rightarrow \mathbb{R}\cup\{+\infty\}$ defined by
$$
\sigma_{\Sigma}(v)=\sup\{\langle\xi,v\rangle;\xi\in \Sigma\}.
$$
\item{$(S_7)$} {\it Let $\Sigma,\Delta$ be nonempty weak$^*$-closed convex subsets of $X^*$. Then
$$
\Sigma \subset \Delta \ \mbox{ if, and only if, } \ \sigma_{\Sigma}(v)\leq
\sigma_{\Delta}(v), \ v\in X;
$$ }
\item{$(S_8)$} {\it $ I^0(u;v)=\max\{\langle\mu,v\rangle; \mu\in \partial I(u)\}$, which means $I^0(u;.)$ is the support function of $\partial I(u)$;}
\item{$(S_9)$} {\it  If $I\in C^1(X,{\mathbb{R}})$, then $\partial (I+J)(u)=\partial I(u)+\partial J(u)$;}
\item{$(S_{10})$} {\it The map $\partial I$ is weak*-closed, in the sense that if \mbox{$\{(u_j,\xi_j)\} \in (X,X^*)$} is a sequence such that $\xi_j\in \partial I (u_j)$ and $u_j \rightarrow u$, then $\xi\in\partial I(u)$;}
\item{$(S_{11})$}~{$(u,v)\mapsto I^0(u;v)$ \it is upper semicontinuous;}
\item{$(S_{12})$}~{\it The set-valued map $x\mapsto\partial f(x)$ is upper semicontinuous, that is, for every $x_0\in X$ and every $\epsilon>0$, there is $\delta=\delta(x_0,\epsilon)>0$ such that for each $x$ with $\|x-x_0\|_{X}<\delta $ and for each $\xi\in \partial f(x)$ there is a $\xi_0\in \partial f(x_0)$ such that $\|\xi-\xi_0\|_{X^*}<\epsilon$}

\end{description}
The first result in this article is a version of Chain Rule theorem. We point out that this result is a different version of that in the books of Rockafellar and Wets \cite{Rokafellar} and Clarke \cite{clarke}, where one has $g$ of class $C^1$ and $f$ a locally Lipschitz function. We observe that we didn't find this kind of result in the literature.

\begin{theo}\label{RegradaCadeia}
Let $f:X \rightarrow \mathbb{R}$, and let $g:\mathbb{R}\to \mathbb{R}$
be function. Suppose that $f$ is Lipschitz in a neighborhood of $x$ and that $g$ is strictly differentiable near $f(u)$. Then $F = g\circ f$ is Lipschitz in a neighborhood of $x$, and one has
\begin{equation}\label{partial(1)}
\partial F(u)\subset  g'(f(u))\partial f(u),\ u\in X.
\end{equation}
\end{theo}
\begin{rem}
The meaning of (\ref{partial(1)}) is that every element
$z$ of $\partial F(u)$ can be represented as
$$
\langle z,v\rangle= g'(f(u))(\langle \xi,v\rangle), \ v\in X,
$$
for some $\xi\in \partial f(u)$.
\end{rem}
{\it Proof of theorem \ref{RegradaCadeia}:} We begin noticing that
\[
g'(f(x)).\partial f(x) = \{ g'(f(x)).\xi \in X^*; \xi \in \partial f(x) \}, \, x \in X.
\]
For every $x, v \in X$, let
\[
q_0(x, v) := \max \{g'(f(x))\langle \xi, v \rangle; \xi \in \partial f(x) \}.
\]
The relation above defines a support function for $g'(f(x)).\partial f(x)$.
\begin{claim}
\label{claim1}
The set $g'(f(x)).\partial f(x) \subset X^*$ is convex and weakly-* closed.
\end{claim}
The convexity of $g'(f(x)).\partial f(x)$ follows from that of $\partial f(x)$. For the \mbox{weak-*} closedness, let $(\hat \xi_n) \subset g'(f(x)).\partial f(x)$ be such that $\hat \xi_n \stackrel{*}{\rightharpoonup} \hat\phi$, where $\stackrel{*}{\rightharpoonup}$ stands for weak* convergence. Note that $\hat \xi _n = g'(f(x)).\xi_n$, where $(\xi_n) \subset \partial f(x)$. Thus $\|\xi_n\|_* \le K(x)$, where $K(x) > 0$ is the Lipschitz constant of $f$ at $x$, and $(\xi_n) \subset \overline{B_{K(x)}}$. Since $\overline{B_{K(x)}}$ is weakly-* compact, there is a subsequence $(\xi_j) \subset (\xi_n)$ such that $\xi_j \stackrel{*}{\rightharpoonup} \xi \in \overline{B_{K(x)}}$. Using the weak-* closedness of $\partial f(x)$, we have $\xi \in \partial f(x)$. Then $\hat \xi_j = g'(f(x))\xi_j \stackrel{*}{\rightharpoonup} g'(f(x))\xi \in g'(f(x))\partial f(x)$, and, by uniqueness of weak-* limit, $\hat \phi = g'(f(x))\xi \in g'(f(x)).\partial f(x)$, finishing the proof of the claim.

Next, let $F^0(x;v)$ and $q_0(x;v)$ be the support functions of $\partial F(x)$ and $g'(f(x)).\partial f(x)$, respectively. Then, to demonstrate the desired inclusion, it suffices to show that $F^0(x;v) \le q_0(x;v)$, for all $x, v \in X$ (see $S_7$). To this end, consider $\epsilon > 0$ and $v \neq 0$ (since for $v = 0$, $F^0(x;0) = 0 = q_0(x;0)$ and the conclusion follows). Now define $q_\epsilon := \sup\{ g'(y).\langle \xi, v \rangle ; \xi \in \partial f(B_\epsilon(x)), y \in (f(x) - \epsilon, f(x) + \epsilon) \}$. Observe that $F \in Lip_{loc}(X, \mathbb R)$, since $g \in C^1$. By definition of $\limsup$, there are $h_\epsilon \in X, \lambda_\epsilon > 0$, with $h_\epsilon \rightarrow 0$ and $\lambda_\epsilon \rightarrow 0^+$, such that $x + h_\epsilon$, $x + h_\epsilon + \lambda_\epsilon v \in B_\epsilon (x)$,
\[
F^0(x;v) - \epsilon \le \frac{F(x + h_\epsilon + \lambda_\epsilon v) - F(x + h_\epsilon)}{\lambda_\epsilon}
\]
and $f(x + h_\epsilon + \lambda_\epsilon v), f(x + h_\epsilon) \in (f(x) - \epsilon, f(x) + \epsilon)$, since $f$ is continuous. By Mean Value theorem, there is $z_\epsilon$ such that
\[
g(f(x + h_\epsilon + \lambda_\epsilon v)) - g(f(x + h_\epsilon)) = g'(z_\epsilon)(f(x + h_\epsilon + \lambda_\epsilon v)-f(x+h_\epsilon)).
\]
On the other hand, by Lebourg's Mean Value theorem,
\[
f(x + h_\epsilon + \lambda_\epsilon v) - f(x + h_\epsilon) = \langle \xi_\epsilon, \lambda_\epsilon v\rangle,
\]
for $\xi_\epsilon \in \partial f(w_\epsilon)$, with $w_\epsilon \in [x + h_\epsilon, x + h_\epsilon + \lambda_\epsilon v]$. 
Then
\begin{equation}\label{desv1}
F^0(x;v) - \epsilon \le g'(z_\epsilon)\langle \xi_\epsilon, v\rangle \le q_\epsilon (x;v),
\end{equation}
since, in particular, $z_\epsilon \in (f(x) - \epsilon, f(x) + \epsilon)$.
\begin{claim}
\label{claim2}
$\lim_{\epsilon \rightarrow 0}q_\epsilon(x;v) = q_0(x;v)$, for $x, v \in X$.
\end{claim}
It suffices to show that given $\sigma' > 0$, there is $\epsilon_0 > 0$ such that
\[
q_0(x, v) - \sigma' \le q_\epsilon (x,v) \le q_0(x, v) + \sigma', \ \forall \epsilon \in (0, \epsilon_0).
\]
Observe that
\begin{equation}
\label{1.6}
q_0(x,v) - \sigma' \le q_0(x,v) \le q_\epsilon (x,v).
\end{equation}
Let us fix $\delta > 0$ and choose $\epsilon_1 > 0$ such that
\begin{equation}
\label{1.7}
f(B_\epsilon(x)) \subset (f(x) - \delta, f(x) + \delta), \ \forall \epsilon \in (0, \epsilon_1).
\end{equation}
By the continuity of $g'$, there is $\rho > 0$ such that
\begin{equation}
\label{1.8}
g'\big( (f(x) - \rho, f(x) + \rho) \big) \subset \big( g'(f(x)) - \delta, g'(f(x)) + \delta \big).
\end{equation}
Since $\partial f$ is lower semicontinuous (see $S_{12}$), there is $\epsilon_2 > 0$ such that for all $y \in B_{\epsilon_2}(x)$ and $ \xi \in \partial f(x)$, there is $\xi_0 \in \partial f(y)$ verifying $\|\xi - \xi _0 \| \le \delta$,
that is,
\begin{equation}
\label{1.9}
\partial f(B_{\epsilon_2}(x)) \subset B_\delta + \partial f(x) = B_\delta(\partial f(x)).
\end{equation}
For $\epsilon \in (0, \min\{\epsilon_1, \epsilon_2, \rho\})$, by \eqref{1.8} and \eqref{1.9},
\begin{align*}
q_\epsilon (x, v) & = \max\{g'(y). \langle \xi, v \rangle ; \xi \in \partial f(B_\epsilon(x)), y \in (f(x) - \epsilon, f(x) + \epsilon)\}\\
    & \le \sup\{g'(y). \langle\xi, v \rangle ; \xi \in B_\delta + \partial f(x), g'(y) \in (g'(f(x)) - \delta, g'(f(x)) + \delta)\}.
\end{align*}
Thus
\begin{align*}
q_\epsilon (x, v) \le & \sup\{g'(y)\langle \xi, v \rangle; \xi \in B_\delta, g'(y) \in (-\delta, \delta)\}\\
                       & + \sup\{g'(y)\langle \xi, v \rangle; \xi \in B_\delta, g'(y) \in \{ g'(f(x)) \} \ \}\\
                       & + \sup\{g'(y)\langle \xi, v \rangle; \xi \in \partial f(x), g'(y) \in (-\delta, \delta)\}\\
                       & + \sup\{g'(y)\langle \xi, v \rangle; \xi \in \partial f(x), g'(y) \in \{ g'(f(x)) \} \ \},
\end{align*}
which implies
\begin{align*}
q_\epsilon (x, v) \le & \sup\{|g'(y)|. \|\xi\|_*\|v\|; \xi \in B_\delta, g'(y) \in (-\delta, \delta)\}\\
                      & + \sup\{|g'(y)|. \|\xi\|_*\|v\|; \xi \in B_\delta, g'(y) \in \{ g'(f(x)) \} \ \}\\
                      & + \sup\{|g'(y)|. \|\xi\|_*\|v\|; \xi \in \partial f(x), g'(y) \in (-\delta, \delta)\}\\
                      & + q_0(x, v)\\
                  \le & \delta^2 \|v\| + |g'(y)|\delta \|v\| + \delta K(x)\|v\| + q_0(x, v).
\end{align*}
Taking the limit $\delta \rightarrow 0^+$,
\begin{equation}
\label{1.10}
q_\epsilon(x,v) \le q_0(x,v) \le q_0(x, v) + \sigma'.
\end{equation}
Therefore we conclude Claim \ref{claim2} by \eqref{1.6} and \eqref{1.10}.

Finally, taking the limit $\epsilon \rightarrow 0^+$ on (\ref{desv1}),
\[
F^0(x, v) \le q_0(x,v),
\]
and, by Claim \ref{claim1}, if follows that
\[
\partial F(x) \subset g'(f(x)).\partial f(x).
\]
\hfill\rule{2mm}{2mm}

As in \cite{Chang}, we say that $I$ satisfies the  {\bf nonsmooth Palais-Smale condition
at level $\boldsymbol{c\in\mathbb{R}}$} ({\bf nonsmooth $\boldsymbol{(PS)_c-}$condition} for
short), if every sequence $\{u_n\}\subset X$ such that
$$
I(u_n)\rightarrow c \ \text{and}\ m_{I}(u_n)\rightarrow 0,
$$
has a strongly convergent subsequence.


\section{Nonsmooth Ricceri's multiplicity theorem and a review on Orlicz-Sobolev spaces}\label{SectionOrlicz}

The next theorem is similar to \cite[Theorem 2.1]{Kristaly}  and its proof will be omitted.
For every $\tau\geq0$, let $G_\tau$ denote the class of functions
$$
\text{G}_{\tau} = \{g\in C^1(\mathbb{R},\mathbb{R}) \text{ is bounded, and } g(t)=t \text{ for any } t\in [-\tau,\tau]\}.
$$

\begin{theo}\label{Theo2.1}
Let $(X,\parallel.\parallel)$ be a real reflexive Banach space and $\widetilde{X}_i$ $i=1,2$ be two Banach spaces such that the embeddings $X\hookrightarrow \widetilde{X}_i$ are compact. Let $\Gamma$ be a real interval, $\Phi_1:X\rightarrow\mathbb{R}_+$ is sequentially weakly lower semicontinuous,  such that
$$
\tilde\eta_0(\parallel u\parallel)\leq \Phi_1(u)\leq \tilde\eta_1(\parallel u\parallel), \ u\in X,
$$
with $\tilde\eta_i:[0,\infty)\rightarrow [0,\infty)$ ($i=1,2$) nondecreasing
and let $\Phi_{i+1}:\widetilde{X}_i\rightarrow\mathbb{R}$ ($i=1,2$) be two locally Lipschitz functions such that
$$
\widehat{J}_{\lambda,\mu}=\Phi_1+\lambda\Phi_2+\mu g\circ\Phi_3
$$
restricted to $X$ satisfies the $(PS)_c$-condition for every $c\in \mathbb{R}$, $\lambda\in \Gamma$, $\mu\in [0,|\lambda|+1]$ and $g\in\mbox{G}_\tau$, $\tau\geq0$. Assume that $\Phi_1+\lambda \Phi_2$ is coercive on $X$ for all $\lambda\in \Gamma$, and that there is $\rho\in\mathbb{R}$ such that
$$
\displaystyle\sup_{\lambda\in~\Gamma}\inf_{x\in~ X}[\Phi_1(u)+\lambda(\Phi_2(u)+\rho)]<\inf_{u\in~X}\sup_{\lambda\in ~\Gamma}[\Phi_1(u)+\lambda(\Phi_2(u)+\rho)].
$$
Then there are a non-empty open set $A\subset \Gamma$ and $r>0$ with the property that for every $\lambda\in A$ there is $\mu_0\in(0,|\lambda|+1]$ such that, for each $\mu\in[0,\mu_0]$, the functional
$$
J_{\lambda,\mu}=\Phi_1+\lambda\Phi_2+\mu\Phi_3
$$
has at least three critical points in $X$ whose norms are less than $r$.
\end{theo}


Let $a$ be a real-valued function defined on $[0,\infty)$ and having
the following properties:

\noindent $(i)$ \  $a(0)=0$, $a(t)>0$ if $t>0$, and
$\displaystyle\lim_{t\rightarrow \infty}a(t)=\infty$;

\noindent $(ii)$ \  $a$ is nondecreasing, that is, $s>t$ implies
$a(s) \geq a(t)$;

\noindent $(iii)$ \ $a$ is right continuous, that is,
$\displaystyle\lim_{s\rightarrow t^{+}}a(s)=a(t)$.

Then the real-valued function $A$ defined on $\mathbb{R}$ by
$$
A(t)= \displaystyle\int^{|t|}_{0}a(s) \ ds
$$
is called an N-function. For a N-function $A$ and  an open set
$\Omega \subseteq \mathbb{R}^{N}$, the Orlicz space
$L_{A}(\Omega)$ is well known (see \cite{Adams}). When $A$
satisfies $\Delta_{2}$-condition, that is, when there are $t_{0}\geq
0$ and $K>0$ such that $A(2t)\leq KA(t)$, for all $t\geq
t_{0}$,  the space $L_{A}(\Omega)$ is the vectorial space of the
measurable functions $u: \Omega \to \mathbb{R}$ such that
$$
\displaystyle\int_{\Omega}A(|u|) \ dx < \infty.
$$
The space $L_{A}(\Omega)$ endowed with Luxemburg norm
$$
|u|_{A}= \inf \biggl\{\lambda >0:
\int_{\Omega}A\Big(\frac{|u|}{\lambda}\Big)\ dx\leq 1\biggl\}
$$
is a Banach space. The complement function of $A$, denoted by
$\widetilde{A}$, is given by its the Legendre transformation
$$
\widetilde{A}(s)=\displaystyle\max_{t \geq 0}\{st -A(t)\} \ \
\mbox{for} \ \ s \geq 0.
$$
One can show that $A$ is the complement of $\widetilde{A}$, and we have
$$
st \leq A(t) + \widetilde{A}(s).
$$
Using the above inequality, known as Young's inequality, it is possible to prove the following H\"{o}lder type inequality
$$
\biggl|\displaystyle\int_{\Omega}u v \ dx \biggl|\leq
2|u|_{A}|v|_{\widetilde{A}},\,\,\, \forall \,\, u \in
L_{A}(\Omega) \,\,\, \mbox{and} \,\,\,  v \in
L_{\widetilde{A}}(\Omega).
$$
It is worth noticing that the $L_{{A}}(\Omega)$ is reflexive if, and only if, $A$ and $\tilde{A}$ satisfy the $\Delta_2$-condition, with
\[
\left( L_{{A}}(\Omega),|\cdot |_{{A},\Omega}\right)^* \!=\! \left( L_{\tilde{A}}(\Omega),\|\cdot \|_{\tilde{A},\Omega}\right)\ , \
 \left( L_{\tilde{A}}(\Omega),|\cdot |_{\tilde{A},\Omega}\right)^* \!=\!
\left( L_{{A}}(\Omega),\|\cdot \|_{{A},\Omega}\right).
\]

The  Orlicz-Sobolev space $ W^{1,A}(\Omega)$, also denoted by $ W^{1}L_A(\Omega)$, is defined in the same way of Sobolev spaces. The usual Orlicz-Sobolev norm of $ W^{1,A}(\Omega)$ is
$$
\|u\|=|u|_A+|\nabla u|_A.
$$
It follows that  $W^{1,A}(\Omega)$ is a reflexive and separable Banach space, since $A$ and $\tilde{A}$ satisfy the $\Delta_2$-condition. Moreover, the $\Delta_2$-condition also implies that
\begin{equation} \label{CV0}
u_n \to u \,\,\, \mbox{in} \,\,\, L_{A}(\Omega) \Longleftrightarrow \int_{\Omega}A(|u_n-u|) \to 0
\end{equation}
and
\begin{equation} \label{CV1}
u_n \to u \,\,\, \mbox{in} \,\,\, W^{1,A}(\Omega) \Longleftrightarrow \int_{\Omega}A(|u_n-u|) \to 0 \,\,\, \mbox{and} \,\,\, \int_{\Omega}A(|\nabla u_n- \nabla u|) \to 0.
\end{equation}

Another important function related to function $\Phi$, is the
Sobolev conjugate function $A_{*}$ of $A$ defined by
$$
A^{-1}_{*}(t)=\displaystyle\int^{t}_{0}\displaystyle\frac{A^{-1}(s)}{s^{(N+1)/N}}ds, \ t>0.
$$

Let $\Omega$ be a smooth bounded domain of $\mathbb{R}^{N}$. If $\Psi$ is any N-function increasing essentially more slowly than $A_{*}$ near infinity, then the imbedding $W^{1}L_{A}(\Omega)\hookrightarrow L_\Psi(\Omega)$ exists and is compact (see \cite{Adams}).

Moreover, we have the following result, due to \cite[Theorem 6.1(i)]{cianchi}.

\begin{prop}\label{cianchi} Let $\Omega$ be a Lipschitz domain in $\mathbb R^N$, $1 \le \kappa < N$ and $A$ be a Young function such that
$$
\int_{0}^{1} \left(\frac{t}{A(t)}\right)^{\frac{\kappa}{N-\kappa}} dt < \infty.
$$
Let
\begin{equation}
\label{*1}
A_T (t) := \int_{0}^{H^{-1}(t)} \left( \frac{A(\tau)}{\tau} \right)^{\frac{N-2}{N-1}} H(\tau)^{\frac{1}{1-N}}d\tau, \ t \ge 0,
\end{equation}
where
\begin{equation}
\label{*2}
H(r) := \left( \int^r_0 \left( \frac{t}{A(t)} \right)^{\frac{1}{N-1}} dt \right)^{\frac{N-1}{N}}.
\end{equation}
Assume that
$$
\int_{1}^\infty \left(\frac{t}{A(t)}\right)^{\frac{\kappa}{N-\kappa}} dt = \infty.
$$
Then there is a constant $C = C(\Omega, \kappa) > 0$ such that
\begin{equation}\label{6.2}
\|\mathrm{Tr}\,u\|_{L^{A_T}(\partial \Omega)} \le C\|u\|_{W^{\kappa, A}(\Omega)},
\end{equation}
for every $u \in W^{\kappa, A}(\Omega)$. Moreover, $L^{A_T}(\partial\Omega)$ is the optimal Orlicz space in \eqref{6.2}, in the sense that if \eqref{6.2} holds with $A_T$ replaced by another Young function $\Upsilon$, then $L^{A_T}(\partial \Omega) \rightarrow L^{\Upsilon}(\partial \Omega)$.
\end{prop}

\section{Application to a differential inclusion with nonhomogeneous boundary condition}
For $\lambda,\mu>0$, we consider the following differential
inclusion problem, with nonhomogeneous Neumann condition:
$$
\left \{ \begin{array}{l}
-div\big(\phi(\mid\nabla u\mid)\nabla u\big)-\phi(|u|)u\in  \lambda\partial F(u(x))  \ \mbox{in}\ \Omega,\\
\frac{\partial u}{\partial \nu}\in \mu\partial G(u(x)) \text{ on } W^{1,\Phi}(\Omega),                                    
\end{array}
\right.\leqno{(P_{\lambda,\mu})}
$$
where $\Omega\subset\mathbb{R}^{N}$ is a bounded smooth domain with
$N \geq 2$ and the function $\phi(t)t$ is increasing in
$(0,+\infty)$, that is,
$$
\leqno{(\phi_1)} \ \            \displaystyle  \lim_{s\rightarrow 0^+}s\phi (s) =0 \,\,\, \mbox{and} \,\,\, \displaystyle  \lim_{s\rightarrow +\infty} s\phi (s) = + \infty,
$$
$$
\leqno{(\phi_2)} \ \ \  s \mapsto s\phi (s)~\mbox{ is increasing in}~ [0, \infty),
$$
$$
\leqno{(\phi_3)} \ \  \   \ell \leq \frac{\phi(t)t^2}{\Phi(t)} \leq m, \   t>0, \ l,m>0,
$$
where $\Phi(t)=\displaystyle\int^{|t|}_{0}\phi(s)s \, ds$, $l\leq m < l^{*}$, $l^{*}=\displaystyle\frac{lN}{N-l}$ and $m^{*}=\displaystyle\frac{mN}{N-m}$. \\

With these hypotheses, we have that $\Phi$ is a $N$-function satisfying the \mbox{$\Delta_{2}$-condition}. Next, we show some examples of functions $\Phi$, for which the related function $\phi(t)t=\Phi'(t)$, for $t \geq 0$, verifies the conditions $(\phi_1)-(\phi_3)$:
$$
\begin{array}{l}
i) \,\, \Phi(t)=|t|^{p}, \,\,\, \mbox{for} \,\,\, 1<p<N;\\
\mbox{}\\
ii) \,\, \Phi(t)=|t|^{p}+|t|^{q}, \,\,\, \mbox{for} \,\,\, 1<p\leq q<N \,\,\, \mbox{and} \,\,\, q \in (p,p^{*}), \,\,\, \mbox{with} \,\,\, \displaystyle p^{*}=\frac{Np}{N-p};\\
\mbox{}\\
iii) \,\, \Phi(t)=(1+|t|^2)^{\gamma}-1, \,\,\, \mbox{for} \,\,\, \gamma \in (1, \min\{\frac{N}{2},\frac{N}{N-2}\} );\\
\mbox{}\\
iv) \,\,  \Phi(t)=|t|^{p}ln(1+|t|), \,\,\, \mbox{for} \,\,\, 1<p_0<p<N-1
\,\,\, \mbox{with} \,\,\, \displaystyle p_0=\frac{-1+\sqrt{1+4N}}{2}.
\end{array}
$$\\
From now on, we will assume that $F:\mathbb{R}\rightarrow\mathbb{R}$ is a locally
Lipschitz function satisfying
\begin{description}
\item{$(F_1)$} there is $c_2>0$ such that
\[
|\xi | \leq c_1(1+b(|t|)|t|), \ \xi \in \partial F(t), \ t\in \mathbb{R},
\]
with $b:(0,+\infty]\rightarrow \mathbb{R}$ a $C^{1}$ function such that
$$
m< b_0\leq \frac{b(t)t^2}{B(t)}\leq b_1 <l^*,\eqno{(b_1)}
$$
for all $t>0$, with
$$
(b(t)t)'>0, \ t>0, \eqno{(b_2)}
$$
and $B(t)=\int_0^{|t|}b(s)s~ds;$\\
\item{$(F_2)$} $\displaystyle \lim_{t\to 0}\frac{\max\{|\xi|;\xi\in \partial F(t)\}}{\phi(|t|)|t|}=0$;

\item{$(F_3)$} $\displaystyle \limsup_{|t|\to+\infty}\frac{F(t)}{\Phi(t)}\leq 0 ;$

\item{$(F_4)$} $F(0)=0$ and there is $t_0\in\mathbb{R}$ such that
$$
F(t_0)>0.
$$
\end{description}
Let also $G:\mathbb{R}\rightarrow\mathbb{R}$ be another locally Lipschitz function satisfying
\begin{description}
\item{$(G_1)$} there is $c_2>0$ such that
\[
|\xi|\leq c_2(1+\overline{b}(|s|)|s|), \ \xi \in \partial G(s), \ s\in \mathbb{R},
\]
where $\overline{b}:(0,+\infty)\rightarrow(0,+\infty)$ verifies $(b_1),(b_2)$, with $\overline{B}(t)=\int_0^{|t|}\overline{b}(s)s~ds,$ $b_0=\overline{b}_0$ and $b_1=\overline{b}_1$, $1 <\overline{b}_0\leq \overline{b}_1<\overline{l}^*=\frac{l(N-1)}{N-l}$.\
\end{description}

Since $F$ and $G$ are locally Lipschitz, by $(F_1)$ and $(F_2)$, \mbox{$\Phi_2:L_{B}(\Omega)\rightarrow\mathbb{R}$} and $\Phi_3:L_{\overline{B}}(\partial\Omega)\rightarrow\mathbb{R}$ given by
$$
\Phi_2(u)=-\int_\Omega F(u)dx, \
\Phi_3(u)=-\int_{\partial\Omega}G(u)dx \
$$
are well-defined locally Lipschitz functionals. Moreover, due to Alves, Gon\c calves and Santos \cite{Abrantes}, we have
\begin{equation}\label{eq_inclusionF}
\partial \Phi_2(u)\subset-\int_\Omega \partial F(u)dx, \ u\in L_B(\Omega)
\end{equation}
and
\begin{equation}\label{eq_inclusionG}
\partial\Phi_3(u)\subset-\int_{\partial\Omega}\partial G(u)dx \ u\in L_{\overline{B}}(\partial\Omega).
\end{equation}
Let $J_{\lambda,\mu}:W^{1,\Phi}(\Omega)\rightarrow\mathbb{R}$  given by
$$
J_{\lambda,\mu}(u)=\Phi_1(u)+\lambda\Phi_2(u)+\mu\Phi_3(u),\ u\in W^{1,\Phi}(\Omega),
$$
be the energy functional associated to the problem $(P_{\lambda,\mu})$, where \mbox{$\Phi_1(u)=\int_\Omega \Phi(|\nabla u|) + \Phi(u)~dx$} is a $C^1$ functional with derivative
$$
\langle \Phi'_1(u),v\rangle=\int_\Omega\phi(\mid\nabla u\mid)\nabla u\nabla v~dx +\int_\Omega\phi(\mid u\mid) u v~dx, \ u,v\in W^{1,\Phi}(\Omega).
$$
Note that $J_{\lambda,\mu}\in Lip_{loc}(W^{1,\Phi}(\Omega),\mathbb{R})$. To prove the above theorem, we need the following auxiliary results. Under assumptions $(\phi_{1})-(\phi_{3})$, the inequalities listed in the following lemmas are valid. For demonstrations, see \cite{Fukagai}.

\begin{lemma}\label{desigualdadeimportantes}
Let, for $t \ge 0$, \mbox{$\xi_{0}(t)=\min\{t^{l},t^{m}\}$} and
\mbox{$\xi_{1}(t)=\max\{t^{l},t^{m}\}$}. Then
\begin{eqnarray*}
&\Phi(t)\xi_0(\rho)\leq\Phi(t\rho )\leq\Phi(t)\xi_1(\rho),&\\
&\xi_{0}(\|u\|_{\Phi})\leq \displaystyle\int_{\Omega}\Phi(|\nabla u|)
\ dx \ \leq\xi_{1}(\|u\|_{\Phi}), \quad \forall u \in W^{1, \Phi}(\Omega).&
\end{eqnarray*}
\end{lemma}

\begin{lemma}\label{desigualdadeimportantes1}
Let, for $t \ge 0$, \mbox{$\eta_{0}(t)=\min\{t^{b_0},t^{b_1}\}$},
\mbox{$\eta_{1}(t)=\max\{t^{b_0},t^{b_1}\}$}, \mbox{$\overline\eta_{0}(t)=\min\{t^{\overline b_0},t^{\overline b_1}\}$} and \mbox{$\overline\eta_{1}(t)=\max\{t^{\overline b_0},t^{\overline b_1}\}$}. Then
\begin{eqnarray*}
&\eta_{0}(|u|_{B}) \leq \displaystyle\int_{\Omega}B(u)dx \leq \eta_{1}(|u|_{B}),&\\
&\overline\eta_{0}(|u|_{\overline B}) \leq \displaystyle\int_{\Omega}\overline B(u)dx \leq \overline\eta_{1}(|u|_{\overline B}),&\\
&B(s)\eta_{0}(t) \leq B(st) \leq B(s)\eta_{1}(t), \ s,t \in \mathbb{R},&\\
&\overline B(s)\overline\eta_{0}(t) \leq \overline B(st) \leq \overline B(s)\overline\eta_{1}(t), \ s,t \in \mathbb{R}.&
\end{eqnarray*}
\end{lemma}

\begin{lemma}\label{DESIGUALD}For every $s>0$, we have
$$
\widetilde{\Phi}\left(\frac{\Phi(s)}{s}\right)\leq \Phi(s) \quad\text{ and }\quad \widetilde{\Phi}(\phi(s)s)\leq\Phi(2s).
$$
\end{lemma}

The next result is a version of Brezis-Lieb's Lemma
\cite{Brezis_Lieb} for Orlicz-Sobolev spaces and the proof can be
found in \cite{Gossez}.

\begin{lemma}\label{brezislieb} Let $\Omega\subset \mathbb{R}^N$ open set and $\Phi:\mathbb{R}\rightarrow [0,\infty)$
an N-function satisfies $\Delta_2-$condition. If the
complementary function $\widetilde{\Phi}$ satisfies
$\Delta_2-$condition, $(f_n)$ is a bounded sequence in $L_\Phi(\Omega)$ such
that
$$
f_n(x)\rightarrow f(x)\  \text{a.e. }x\in\Omega,
$$
then
$$
f_n\rightharpoonup f \ \text{in }L_\Phi(\Omega).
$$
\end{lemma}

\begin{coro}\label{imersao}
The embeddings $W^{1,\Phi}(\Omega)\hookrightarrow L_{B}(\Omega)$ and 
$W^{1,\Phi}(\Omega)\hookrightarrow L_{\Phi}(\Omega)$ are compact.
\end{coro}
\text{\it Proof:} It is sufficient to show that $B$ increases essentially
more slowly than $\Phi_{*}$ near infinity. Indeed,
$$
\frac{B(kt)}{\Phi_*(t)} \leq \frac{B(1)\eta_1(kt)}{\xi_2(t)} = B(1)k^{b_1}t^{b_1-l^*}, \ k > 0.
$$
Since $b_1<l^*$, we get
$$
\displaystyle\lim_{t\rightarrow +\infty}\frac{B(kt)}{\Phi_*(t)}= 0.
$$
Analogously, we get the same conclusion for $\Phi$.
\hfill\rule{2mm}{2mm}

\begin{theo}\label{imersao1}
The embedding $W^{1,\Phi}(\Omega)\hookrightarrow L_{\overline{B}}(\partial\Omega)$ is compact.
\end{theo}
\text{\it Proof:} Following Theorem \ref{cianchi}, we need to show that $\overline B$ increases essentially more slowly than $\Phi_T$ near infinity.

\begin{claim}\label{phit} There are constants $\overline c_1, \overline c_2 \in \mathbb{R}$, $\overline c_2 > 0$, such that, for $t$ large enough,
\begin{equation}\label{estimativa}
\Phi_T(t) \ge \overline c_1 + \overline c_2 t^{\overline l^*}.
\end{equation}
\end{claim}

Assuming the claim is proved, we have, by Lema \ref{desigualdadeimportantes1}, for $t$ large enough,
\begin{align*}
\frac{\overline B(kt)}{\Phi_T(t)} \le \frac{\overline B(1) \overline\eta_1(kt)}{\overline c_1 + \overline c_2 t^{\overline l^*}} \le \frac{\overline B(1)k^{\overline b_1}t^{\overline b_1}}{\overline c_1 + \overline c_2 t^{\overline l^*}},
\end{align*}
which goes to zero as $t$ goes to infinity. Then $\overline B$ increases essentially more slowly than $\Phi_T$ near infinity and, by Theorem \ref{cianchi}, we have $W^{1,\Phi}(\Omega)\hookrightarrow L_{\overline{B}}(\partial\Omega)$ compactly.

Now we prove the claim. By Lemma \ref{desigualdadeimportantes},
$$
\Phi(1) t^l \le \Phi(t) \le \Phi(1)t^m, \ t \ge 1,
$$
and then
\begin{equation}
\label{*3}
\frac{t^{1-m}}{\Phi(1)} \le \frac{t}{\Phi(t)} \le \frac{t^{1-l}}{\Phi(1)}, \ t \ge 1.
\end{equation}
By \eqref{*2} and \eqref{*3},
\begin{equation}
\label{*4}
\left(\int_0^r \left(\frac{t^{1-m}}{\Phi(1)}\right)^{\frac{1}{N-1}}\right)^{\frac{N-1}{N}} \le H(r) \le \left(\int_0^r \left(\frac{t^{1-l}}{\Phi(1)}\right)^{\frac{1}{N-1}}\right)^{\frac{N-1}{N}}, \ r \ge 1.
\end{equation}
Note that
\begin{align}
\left( \int^r_0 \left( \frac{t^{1-l}}{\Phi(1)} \right)^{\frac{1}{N-1}}dt \right)^{\frac{N-1}{N}} & = \left( \frac{1}{\Phi(1)^{\frac{1}{N-1}}} \right)^{\frac{N-1}{N}} \left( \int^r_0 t^{\frac{1-l}{N-1}} dt\right)^{\frac{N-1}{N}} \nonumber \\
                        & = \frac{1}{\Phi(1)^{\frac{1}{N}}}\left( \frac{1}{1 + \frac{1-l}{N-1}} r^{\frac{N-l}{N-1}}\right)^{\frac{N-1}{N}} \nonumber \\
                        & = \frac{1}{\Phi(1)^{\frac{1}{N}}}\left( \frac{N-1}{N-l} \right)^{\frac{N-1}{N}}r^{\frac{N-l}{N}}. \label{*5}
\end{align}
Analogously, we have
\begin{equation}
\label{*6}
\left( \int^r_0 \left( \frac{t^{1-m}}{\Phi(1)} \right)^{\frac{1}{N-1}}dt \right)^{\frac{N-1}{N}} = \frac{1}{\Phi(1)^{\frac{1}{N}}}\left( \frac{N-1}{N-m} \right)^{\frac{N-1}{N}}r^{\frac{N-m}{N}}.
\end{equation}
By \eqref{*4} - \eqref{*6},
\begin{equation*}
\frac{1}{\Phi(1)^{\frac{1}{N}}}\left( \frac{N-1}{N-m} \right)^{\frac{N-1}{n}}r^{\frac{N-m}{N}} \le H(r) \le \frac{1}{\Phi(1)^{\frac{1}{N}}}\left( \frac{N-1}{N-l} \right)^{\frac{N-1}{N}}r^{\frac{N-l}{N}}.
\end{equation*}
Letting
\[
c_{m, N} := \frac{1}{\Phi(1)^{\frac{1}{N}}}\left( \frac{N-1}{N-m} \right)^{\frac{N-1}{N}}, \ c_{l, N} := \frac{1}{\Phi(1)^{\frac{1}{N}}}\left( \frac{N-1}{N-l} \right)^{\frac{N-1}{N}},
\]
we can rewrite the last inequalities as
\begin{equation}
\label{*7}
c_{m, N} r^{\frac{N-m}{N}} \le H(r) \le c_{l, N}r^{\frac{N-l}{N}}.
\end{equation}
Setting $c_\Phi := \int^1_0 \left( \frac{\Phi(\tau)}{\tau} \right)^{\frac{N-2}{N-1}} H(\tau)^{\frac{1}{1-N}} d\tau$, $\tilde c_{l, N}^1 := {c_{l, N}}^{\frac{1}{1-N}} \Phi(1)^{\frac{N-2}{N-1}}$, we have, by \eqref{*1} and \eqref{*7}, for $t$ large enough, using the inequality preceding \eqref{*3},
\begin{align*}
\Phi_T(t) & \ge \int^1_0 \left( \frac{\Phi(\tau)}{\tau} \right)^{\frac{N-2}{N-1}} H(\tau)^{\frac{1}{1-N}} d\tau + \int^{H^{-1}(t)}_1 \left( \frac{\Phi(\tau)}{\tau} \right)^{\frac{N-2}{N-1}} H(\tau)^{\frac{1}{1-N}} d\tau\\
       & \ge c_\Phi + \int^{H^{-1}(t)}_1 \Phi(\tau)^{\frac{N-2}{N-1}}{c_{l, N}}^{\frac{1}{1-N}}\tau^{\left\{\frac{l-N}{N(1 - N)} - \frac{N - 2}{N - 1}\right\}}d\tau\\
       & = c_\Phi + {c_{l, N}}^{\frac{1}{1-N}} \int^{H^{-1}(t)}_1 \Phi(\tau)^{\frac{N-2}{N-1}}\tau^{\frac{l - N(N - 1)}{N(N - 1)}}d\tau\\
       & \ge c_\Phi + {c_{l, N}}^{\frac{1}{1-N}} \Phi(1)^{\frac{N-2}{N-1}} \int^{H^{-1}(t)}_1 \tau^{\frac{l(N-2)}{N-1}}\tau^{\frac{l - N(N - 1)}{N(N - 1)}}d\tau\\
       & = c_\Phi + \tilde c_{l, N}^1 \int^{H^{-1}(t)}_1 \tau^{\frac{l + Nl(N-2)}{N(N - 1)} - 1}d\tau\\
       & = c_\Phi + \tilde c_{l, N}^1\left( \frac{1}{\frac{l + Nl(N-2)}{N(N-1)}} \tau^{\frac{l + Nl(N-2)}{N(N-1)}} \right)\Big|^{\tau = H^{-1}(t)}_{\tau = 1}\\
       & = c_\Phi + \tilde c_{l, N}^1 \frac{N(N-1)}{l + Nl(N-2)} H^{-1}(t)^{\frac{l + Nl(N-2)}{N(N-1)}} - \tilde c_{l, N}^1 \frac{N(N-1)}{m + Nl(N-2)}.
\end{align*}
Now, let
\[
\tilde c_{l, N}^2 := c_\Phi - \tilde c_{l, N}^3
\]
and
\begin{align*}
\tilde c^3_{l, N} & := \tilde c^1_{l,N} \frac{N(N-1)}{m + Nl(N-2)}\\
                  & = {c_{l, N}}^{\frac{1}{1-N}} \Phi(1)^{\frac{N-2}{N-1}} \frac{N(N-1)}{l + Nl(N-2)}\\
                  & = \Phi(1)^{\frac{1}{N(N-1)}} \left( \frac{N-1}{N-l} \right)^{-\frac{1}{N}} \Phi(1)^{\frac{N-2}{N-1}} \frac{N(N-1)}{l + Nl(N-2)}\\
                  & = \Phi(1)^{\frac{N(N-2)+1}{N(N-1)}}
                  \frac{N(N-1)^{\frac{N-1}{N}}(N-l)^{\frac{1}{N}}}{l + Nl(N-2)}>0.
\end{align*}
So, for $t$ large enough,
\begin{equation}
\label{*8}
\Phi_T(t) \ge \tilde c^2_{l, N} + \tilde c^3_{l,N}H^{-1}(t)^{\frac{l+Nl(N-2)}{N(N-1)}}.
\end{equation}
On the other hand, by \eqref{*7}, we have, for $t = H(r)$ large enough,
\[
t \le \frac{(N-1)^{\frac{N-1}{N}}}{A(1)^{\frac{1}{N}} (N - l)^{\frac{N-1}{N}}} H^{-1}(t)^{\frac{N-l}{N}},
\]
and then
\begin{equation}
\label{*9}
H^{-1}(t) \ge \frac{\Phi(1)^{\frac{1}{N-l}} (N - l)^{\frac{N-1}{N-l}}}{(N-1)^{\frac{N-1}{N-l}}}t^{\frac{N}{N-l}}.
\end{equation}
From \eqref{*8} and \eqref{*9}, for $t$ large enough,
\begin{equation}
\label{*10}
\Phi_T(t) \ge \tilde c^2_{l,N} + \tilde c^4_{l, N} t^{\frac{l+Nl(N-2)}{(N-l)(N-1)}},
\end{equation}
where
\[
\tilde c^4_{l, N} := \tilde c^3_{l,N} \left( \Phi(1)^{\frac{1}{N-l}} \left(\frac{N-l}{N-1}\right)^{\frac{N-1}{N-l}} \right)^{\frac{l+Nl(N-2)}{(N-l)(N-1)}}.
\]
Now observe that
\begin{align}
\frac{l+Nl(N-2)}{(N-l)(N-1)} & = \frac{-l (N-1) + Nl(N-1)}{(N-l)(N-1)} \nonumber\\
                             &= \frac{l(N-1)}{N-l} =: \overline l ^*. \label{*11}
\end{align}
Finally, by \eqref{*10} and \eqref{*11},
\[
\Phi_T(t) \ge \tilde c^2_{l,N} + \tilde c^4_{l,N} t^{\overline l ^*},
\]
for $t$ large enough.\hfill\rule{2mm}{2mm}

\begin{rem} By Theorem \ref{imersao}, if $\overline B(t)=|t|^p$, with $1\leq p<\overline l^*$, the embedding $W^{1,\Phi} (\Omega)\hookrightarrow L^p(\partial\Omega)$ is compact.
\end{rem}


\begin{prop}\label{prop31} The following limit holds:
$$
\displaystyle \lim_{t\to 0^+}\frac{\inf\{\Phi_2(u);u\in W^{1,\Phi}(\Omega), \ \Phi_1(u)<t\}}{t}=0.
$$
\end{prop}
\text{\it Proof:} Applying Lebourg's mean value theorem and using $(F_1)-(F_2)$, for $\epsilon>0$, there is $K=K_\epsilon>0$ such that
\begin{equation}\label{eq33}
\mid F(t)\mid\leq m\epsilon \Phi(t)+KB(t),\ t\in\mathbb{R}.
\end{equation}
By (\ref{eq33}) and Corollary \ref{imersao}, we have
\begin{equation}\label{eq34}
\Phi_2(u)\geq -\epsilon\ \! m\Phi_1(u)-c_\epsilon \eta_1(\| u\|).
\end{equation}
For $t>0$, set $S_t=\{u\in W^{1,\Phi}(\Omega);\Phi_1(u)<t\}$. Using  (\ref{eq34}), we obtain
$$
0\geq\frac{\Phi_2(u)}{t}\geq -\epsilon\ \! m-\frac{c_\epsilon}{t}t^{\frac{b_0-m}{m}}.
$$
Taking the limit $t\to 0^+$, since $\epsilon>0$ is arbitrary, the desired limit follows.
\hfill\rule{2mm}{2mm}\\

\noindent{\it \textbf{Proof of Theorem \ref{theorem1}:}} For $t > 0$, let
$$
\beta(t)=\inf\left\{\Phi_2(u);u\in W^{1,\Phi}(\Omega),\ \Phi_1(u)<t\right\}.
$$
Note that $\beta(t) \le 0$ and, by Proposition \ref{prop31},
\begin{equation}\label{eq35}
\lim_{t\to 0^+}\frac{\beta(t)}{t}=0.
\end{equation}
Now consider the function $u_0 \equiv t_0$ in $\Omega$,  $t_0$ given by $(F_4)$. Then $u_0\in W^{1,\Phi}(\Omega)$. Note that $t_0\neq 0$ (since $F(0)=0$), so $\Phi_2(u_0)<0$. Therefore it is possible to choose a number $\eta>0$ such that
$$
0<\eta<-\frac{\Phi_2(u_0)}{\Phi_1(u_0)}.
$$
By (\ref{eq35}), there is $t_1\in(0,\Phi_1(u_0))$ such that $-\beta(t_1)<\eta t_1$. Thus,
\begin{equation}\label{eq36}
-\frac{\beta(t_1)}{t_1}<-\frac{\Phi_2(u_0)}{\Phi_1(u_0)}.
\end{equation}
Due to the choice of $t_1$ and using (\ref{eq36}), we conclude that there is $\rho_0>0$ such that
\begin{equation}\label{eq37}
-\beta(t_1)<\rho_0<-\frac{t_1\Phi_2(u_0)}{\Phi_1(u_0)}<-\Phi_2(u_0).
\end{equation}
Now let $\varphi:W^{1,\Phi}(\Omega)\times I\rightarrow\mathbb{R}$ be defined by
$$
\varphi(u,\lambda)=\Phi_1(u)+\lambda\Phi_2(u)+\lambda\rho_0, \ u\in W^{1,\Phi}(\Omega) , \lambda \in I,
$$
where $I=[0,+\infty)$. We assert that $\varphi$ satisfies the inequality
\begin{equation}\label{eq38}
\displaystyle\sup_{\lambda\in I}\inf_{u\in W^{1,\Phi}(\Omega)}\varphi(u,\lambda)<\inf_{u\in W^{1,\Phi}(\Omega)}\sup_{\lambda\in I}\varphi(u,\lambda).
\end{equation}
Indeed, the function
$$
I\ni\lambda\mapsto\displaystyle \inf_{u\in W^{1,\Phi}(\Omega)}\varphi(u,\lambda)
$$
is upper semicontinuous on $I$. It follows from (\ref{eq37}) that
$$
\displaystyle\lim_{\lambda\to+\infty}\inf_{u\in W^{1,\Phi}(\Omega)} \varphi(u,\lambda)\leq \lim_{\lambda\to+\infty}\left(\Phi_1(u_0)+\lambda(\rho_0+\Phi_2(u_0)\right)=-\infty.
$$
Thus we find an element $\lambda_0\in I$ such that
\begin{equation}\label{eq39}
\sup_{\lambda\in I}\inf_{u\in W^{1,\Phi}(\Omega)}\varphi(u,\lambda)=\inf_{u\in W^{1,\Phi}(\Omega)}\varphi(u,\lambda_0).
\end{equation}
Since $-\beta(t_1)<\rho_0$, it follows from the definition of $\beta$ that, for all $u\in W^{1,\Phi}(\Omega)$ with $\Phi_1(u)<t_1$, $-\Phi_2(u)<\rho_0$. Hence
\begin{equation}\label{eq310}
t_1\leq \inf\left\{\Phi_1(u);u\in W^{1,\Phi}(\Omega), \ -\Phi_2(u)\geq \rho_0\right\}.
\end{equation}
On the other hand
\begin{eqnarray}
\inf_{u\in W^{1,\Phi}(\Omega)} \sup_{\lambda\in I}\varphi(u,\lambda) & = & \inf_{u\in W^{1,\Phi}(\Omega)}\left(\Phi_1(u)+\sup_{\lambda\in I}(\lambda\rho_0+\lambda\Phi_2(u))\right)\nonumber\\
                            & = &\inf_{u\in W^{1,\Phi}(\Omega)}\left\{\Phi_1(u);-\Phi_2(u)\geq \rho_0\right\}.\nonumber
\end{eqnarray}
Thus, (\ref{eq310}) is equivalent to
\begin{equation}\label{eq311}
t_1\leq \inf_{u\in W^{1,\Phi}(\Omega)}\sup_{\lambda\in I}\varphi(u,\lambda).
\end{equation}
There are two cases to consider. For $\lambda_0\in [0,\frac{t_1}{\rho_0})$,
$$
\inf_{u\in W^{1,\Phi}(\Omega)}\varphi(u,\lambda_0)\leq\varphi(0,\lambda_0)=\lambda_0\rho_0<t_1.
$$
Combining this inequality with (\ref{eq39}) and (\ref{eq311}), (\ref{eq38}) follows; if $\frac{t_1}{\rho_0}\leq \lambda_0$, then, from (\ref{eq37}),
\begin{eqnarray}
\inf_{u\in W^{1,\Phi}(\Omega)}\varphi(u,\lambda_0)&\leq& \Phi_1(u_0)+\lambda_0(\rho_0+\Phi_2(u_0))\nonumber\\
&\leq& -\frac{t_1}{\rho_0}\Phi_2(u_0)+\frac{t_1}{\rho_0}(\rho_0+\Phi_2(u_0))=t_1 , \nonumber
\end{eqnarray}
and we conclude (\ref{eq38}) by another application of (\ref{eq39}) and (\ref{eq311}). Fix $g\in G_{\tau}$. Now we are in the position to apply Theorem \ref{Theo2.1}: we choose $X=W^{1,\Phi}(\Omega)$, $\widetilde{X}_1=L_B(\Omega)$, $\widetilde{X}_2=L_{\overline{B}}(\partial \Omega)$, $\Gamma=I=[0,+\infty)$, and $c\in \mathbb{R}$. We shall prove that the functional $\widetilde{J}_{\lambda,\mu}:W^{1,\Phi}(\Omega)\rightarrow\mathbb{R}$ given by
$$
\widetilde{J}_{\lambda,\mu}(u)=\Phi_1(u)+\lambda\Phi_2(u)+\mu (g\circ\Phi_3)(u), \ u\in
W^{1,\Phi}(\Omega),
$$
satisfies the $(PS)_c$ condition. By $(S_3)$, $(S_7)-(S_9)$ and Theorem \ref{RegradaCadeia}, we have, for every $u,v\in W^{1,\Phi}(\Omega)$ and $w\in \partial \widetilde{J}_{\lambda,\mu}(u)$,
\begin{equation}\label{eq312}
\langle w,v\rangle = \langle\Phi_1'(u),v\rangle+\lambda \langle\xi_F,v\rangle+\mu g'\left(\Phi_3(u)\right) \langle\xi_G,v\rangle,
\end{equation}
for some $\xi_F\in \partial \Phi_2(u)$ and $\xi_G\in\partial \Phi_3(u)$. First of all, let us observe that $\Phi_1+\lambda\Phi_2$ is coercive on $W^{1,\Phi}(\Omega)$, due to $(F_3)$; thus, the functional $\widetilde{J}_{\lambda,\mu}$ is also coercive on $W^{1,\Phi}(\Omega)$. Consequently, it is enough to consider a bounded sequence $\{u_n\}\subset W^{1,\Phi}(\Omega)$ such that
\begin{equation}\label{eq313}
\widetilde{J}_{\lambda,\mu}(u_n)\to c_{\lambda,\mu}\ \text{and }m_{\lambda,\mu}(u_n)\to 0.
\end{equation}
Assuming $\| w_n \|_*=m_{\lambda,\mu}(u_n)$, there are $\xi^F_n\in\partial \Phi_2(u_n)$ and $\xi_n^G\in\partial \Phi_3(u)$ such that
$$
\langle w_n,v\rangle=\langle \Phi'_1(u_n),v\rangle+\lambda\langle\xi_n^F,v \rangle+\mu g'\left( \Phi_3(u_n)\right)\langle \xi_n^G,v\rangle,
$$
for any $v\in W^{1,\Phi}(\Omega)$. By the boundedness of the sequence $\{u_n\}$ in $W^{1,{\Phi}}(\Omega)$, which is a reflexive space, we have, from Corollary \ref{imersao} and Theorem \ref{imersao1}, $u\in W^{1,\Phi}(\Omega)$ such that
$$
u_n \rightharpoonup u \ \mbox{em}\ W^{1,\Phi}(\Omega),
$$
$$
u_n\to u \ \mbox{em}\ L_{B}(\Omega),
$$
$$
u_n\to u \ \mbox{em}\ L_{\overline{B}}(\partial\Omega)
$$
and
$$
u_n(x) \to u(x) \,\mbox{a.e.} \, x\in\Omega
$$
{\bf Claim 1:} $\{\xi^F_n\}\subset\partial F(u_n)$ is bounded in $L_{\widetilde{B}}(\Omega)$ and $\{g'\left(\Phi_3(u_n)\right)\xi_n^G\}$ is bounded in $L_{\widetilde{\overline{B}}}(\partial \Omega)$, where $\{\xi^G_n\}\subset \partial G(u_n)$.

By $(F_1)$ and Lemmas \ref{desigualdadeimportantes1} and \ref{DESIGUALD}, we have
\begin{eqnarray}\label{eqbounded}
\int_\Omega \widetilde{B}(\xi_n^F)&\leq &c\mid\Omega\mid+c\int_\Omega \widetilde{B}(b(| u_n|)u_n)\nonumber\\
&\leq& c+c\int_\Omega B(u_n)\nonumber\\
&\leq& c+c\eta_1(|u_n|_{B})\nonumber\\
&\leq& c+c\eta(\|u_n\|)\leq c, \ n\in\mathbb{N},
\end{eqnarray}
where $c>0$. Another application of Lemmas \ref{desigualdadeimportantes1} and \ref{DESIGUALD} and using $(G_1)$, H\"{o}lder's inequality and Corollary \ref{imersao} gives
\begin{eqnarray}
\left|\langle \xi^G_n, v\rangle\right| & = & \left|\int_{\partial \Omega}\xi^G_nv\right|\nonumber\\
&\leq&c_2\int_{\partial \Omega}(1+\overline{b}(|u_n|)|u_n|)|v|\nonumber\\
&\leq&c(1+ |u_n|_{\overline{B},\partial
\Omega})|v|_{\overline{B},\partial \Omega}\nonumber\\
&\leq &c(1 + \|u_n\|) \| v \|\leq c\|v\|.\label{b1}
\end{eqnarray}
By Lebourg's Theorem and $(G_1)$, there is $c>0$ such that
\[
|G(t)| \le c(1 + \overline b(|t|)|t|^2).
\]
Thus, arguing as in \eqref{b1},
$$
\left|\Phi_3(u_n)\right|\leq 
c+c\overline\eta_1(|u_n|_{\overline B,\partial \Omega})\leq c+c\overline \eta_1(\|u_n\|)\leq c,
$$
for any $n\in\mathbb{N}$. Since $g\in C^1$, we have
$$
\left|g'(\Phi_3(u_n))\langle\xi_n^G,v\rangle\right| \leq c|g'(\Phi_3(u_n))|~\| v\| \leq c\| v\|,
$$
for any $n\in\mathbb{N}$. Then $\{g'(\Phi_3(u_n))\xi_n^G\}$ is bounded in $L_{\widetilde{\overline{B}}}(\partial\Omega)$.

Using (\ref{eq312}) with $v=u_n-u$, we obtain
\begin{equation}\label{eqon}
o_n(1)=\langle\Phi_1'(u_n),u_n-u\rangle-\lambda\int_\Omega\xi_n^F(u_n-u)-\mu g'(\Phi_3(u_n))\int_{\partial \Omega}\xi_n^G(u_n-u).
\end{equation}
From Claim 1,
\begin{align}
\label{eqinxifg}
\left|\lambda\int_\Omega \xi_n^F (u_n- \right.& u)\left. + \mu g'(\Phi_3(u_n))\int_{\partial \Omega}\xi_n^G(u_n-u)\right| \le \nonumber\\
\le & \lambda |\xi_n^F|_{\widetilde{B}} |u_n-u|_B + \mu \left|g'(\Phi_3(u_n))\xi_n^{G}\right|_{\widetilde{\overline B},\partial \Omega} |u_n-u|_{\overline B,\partial \Omega}\nonumber\\
\le & c\lambda |u_n-u|_{B} + c\mu|u_n-u|_{\overline B,\partial \Omega}\to 0.
\end{align}
By lemma \ref{DESIGUALD}
$$
\int_\Omega\widetilde{\Phi}\left(\phi(u_n)u_n\right)\leq \int_{\Omega}\Phi(2u_n)\leq 2^m\int_\Omega\Phi(u_n)\leq 2^m\xi_1(|u_n|_{\Phi})\leq c,
$$
which implies
\begin{equation}\label{eqconv01}
\left|\int_\Omega \phi(|u_n|)u_n(u_n-u)\right|\leq |\phi(u_n)u_n|_{\widetilde{\Phi}}|u_n-u|_\Phi\leq c|u_n-u|_{\Phi}\to 0.
\end{equation}
Combining (\ref{eqon})-(\ref{eqconv01}),
\begin{equation}\label{eqconv1}
\left|\int_\Omega \phi(\mid\nabla u_n\mid)\nabla u_n\nabla(u_n-u)\right|\to 0.
\end{equation}
Since, for some subsequence,
$$
\frac{\partial u_n}{\partial x_i}
\rightharpoonup \frac{\partial u}{\partial x_i} \text{ in }L_\Phi(\Omega),
$$
we have that
\begin{equation}\label{eqconv2}
\int_\Omega\phi(\mid\nabla u\mid)\nabla u\nabla(u_n-u)\to 0.
\end{equation}
From (\ref{eqconv1}) and (\ref{eqconv2}),
$$
0\leq\int_\Omega\left(\phi(|\nabla u_n|)\nabla u_n-\phi(|\nabla u|)\nabla u\right)\left(\nabla u_n-\nabla u\right)\to 0.
$$
Applying a result due to Dal Maso and Murat \cite{Maso}, we obtain
$$
\nabla u_n\to \nabla u \mbox{ a.e. }\Omega.
$$
Now, using Lebesgue's Theorem, we obtain
\begin{equation} \label{E2}
u_{n} \to u \mbox{ in } W^{1,\Phi}(\Omega).
\end{equation}
It remains to apply Theorem \ref{Theo2.1} in order to obtain the conclusion. \hfill\rule{2mm}{2mm}\\

\vspace{1 cm}

\noindent Rodrigo C. M. Nemer and Jefferson A. Santos \\
\noindent Universidade Federal de Campina Grande,\\
Unidade Acad\^emica de Matem\'atica,\\
\noindent CEP:58109-970, Campina Grande - PB, Brazil\\
\noindent e-mail: rodrigo@mat.ufcg.edu.br and jefferson@mat.ufcg.edu.br \\


\begin{thebibliography}{99}

\bibitem{Adams} A. Adams and J. F. Fournier, {\it Sobolev spaces},  2nd ed., Academic
Press, (2003).

\bibitem{Abrantes} {C.O. Alves, J.V. Gon\c calves and J.A. Santos,} {\it Strongly Nonlinear Multivalued Elliptic
Equations on a Bounded Domain}, J Glob Optim, v.58, p.565-513 (2014).

\bibitem{Rubia} C.O. Alves and R.G. Nascimento, {\it Nonlinear perturbations of a periodic elliptic problem with discontinuousnonlinearity in $\mathbb{R}^N$.} Z. Angew. Math. Phys. 63, 107-124 (2012).

\bibitem{Arcoya} D. Arcoya and J. Carmona, {\it A nondifferentiable extension of a theorem of Pucci-Serrin and applications.} J. Differ. Equ. 235(2), 683-700 (2007).

\bibitem{Bn1} G. Bonanno, {\it Some remarks on a three critical points theorem}. Nonlinear Anal. 54, 651-665 (2003)

\bibitem{Bn2} G. Bonanno, {\it A critical points theorem and nonlinear differential problems.} J. Global Optim. 28(3-4), 249-258 (2004)

\bibitem{Bn3} G. Bonanno and P. Candito, {\it Non-differentiable functionals and applications to elliptic problems with discontinuous nonlinearity.} J. Differ. Equ. 244(12), 3031-3059 (2008)



\bibitem{Brezis_Lieb} H. Brezis. and E. Lieb, {\it A relation between pointwise convergence of functions and
convergence of functinals}, Proc. Amer. Math. Soc. 88,
486-490 (1983).

\bibitem{carl} S. Carl, V. K. Le and D. Motreanu, {\it Nonsmooth variational problems and their inequalities.
Comparison principles and applications}, Springer Monographs in
Mathematics. Springer, New York, (2007).

\bibitem{clarke} F.H. Clarke, {\it Optimization and Nonsmooth Analysis}, John Wiley \& Sons, N.Y, 1983.

\bibitem{Chang} K. C. Chang, {\it Variational methods for
    nondifferentiable functionals and their applications to partial
    differential equations.} J. Math. Analysis Aplic., 80,
    102-129 (1981).
		
\bibitem{cianchi}  A. Cianchi,{\it On some aspects of the Theory of Orlicz-Sobolev Spaces}.  In: Around the Research of Vladimir Maz’ya. I, Volume 11 of International Mathematics Series, NY, pp. 81-104. Springer, New York (2010).

		
\bibitem{dacorogna}  Dacorogna, B., {\it Introduction to the Calculus of Variotions},  Imperial College Press (2009).


\bibitem{Fukagai} N. Fukagai, M. Ito and K. Narukawa,
\emph{Positive solutons of quasilinear elliptic equations with
critical Orlicz-Sobolev nonlinearity on $\mathbb{R}^{N}$},
Funkciallaj Ekvacioj, 49(2006)235-267.

\bibitem{narukawa-2} N. Fukagai and M. Narukawa, {\it Nonlinear eigenvalue problem for a model equation of an elastic surface}, Hiroshima Math. J. 25, 19-41 (1995).

\bibitem{gasinski} L. Gasi\' nski and N.S. Papageorgiou, {\it Nonsmooth Critical Point Theory and Nonlinear Boundary Value Problems.} Chapman \& Hall/CRC, London, 2004.

\bibitem{Gossez} J.P. Gossez {\it Orlicz-Sobolev spaces and nonlinear elliptic boundary
value problems.} In: Fuk, Svatopluk and Kufner, Alois (eds.):
Nonlinear Analysis, Function Spaces and Applications, Proceedings of
a Spring School held in Horn Bradlo, 1978. [Vol 1]. BSB B. G.
Teubner Verlagsgesellschaft, Leipzig, 1979. pp. 59-94.

\bibitem{Kristaly} A. Krist\'aly, W. Marzantowicz and C. Varga, {\it A non-smooth three critical points theorem with applications in
      differential inclusions}, J. Glob. Optim. Volume 46, Number 1 (2010), 49-62.
			
\bibitem{Kristaly2} A. Krist\'aly, {\it Infinitely many solutions for a differential inclusion problem in $\mathbb R ^N$.} J.Differ.Equ. 220, 511-530, 2006.

\bibitem{Marano} S.A. Marano and D. Motreanu, {\it Infinitely many critical points of non-differentiable functions and applications to a Neumann-type problem involving the p-Laplacian.} J. Differential Equations 182, 108-120, 2002.

\bibitem{Maso} G. Dal Maso and F. Murat, {\it Almost everywhere convergence of gradients of solutions to nonlinear elliptic
    systems}, Nonlinear Anal. 31 (1998), 405-412.

\bibitem{Marano2} S.A. Marano and D. Motreanu, {\it On a three critical points theorem for non-differentiable functions and applications to nonlinear boundary value problems.} Nonlinear Anal. 48, 37-52 (2002).


\bibitem{papa-1} D. Motreanu and P.D. Panagiotopoulos, {\it Minimax Theorems and Qualitative Properties of the Solutions of Hemivariational Inequalities}, Nonconvex Optim. Appl. 29 Kluwer, Dordrecht, 1998.


\bibitem{Pucci} P. Pucci and J. Serrin, {\it A mountain pass theorem.} J. Differ. Equ. 60, 142-149 (1985).

 \bibitem{Rao1}  M.N. Rao and Z.D. Ren, {\it Theory of Orlicz Spaces},~ Marcel Dekker, New York, (1985).	

\bibitem{Ricceri} B. Ricceri, {\it Minimax theorems for limits of parametrized functions having at most one local minimum lying in a certain set}.
      Topol. Appl. 153, 3308-3312 (2006).

\bibitem{Ricceri2} B. Ricceri, {\it Infinitely many solutions of the Neumann problem for elliptic equations involving the p-Laplacian.} Bull. London Math. Soc. 33, 331-340, 2001

\bibitem{Rokafellar} R.T. Rockafellar and R.J.-B. Wets, {\it Variational Analysis}. Springer-Verlag, Berlin (1998).


\end{thebibliography}
\end{document}